\documentclass{elsart3-1}

\usepackage{amssymb}

\usepackage[english,francais]{babel}

\newtheorem{e-proposition}[theorem]{Proposition}

\newtheorem{e-definition}[theorem]{Definition\rm}

\renewcommand{\theequation}{\arabic{equation}}
\newcommand{\R}{\mbox{$\mathbb R$}}

\newcommand{\C}{\mbox{$\mathbb C$}}
\newcommand{\Sp}{\mbox{$\mathbb S$}}
\newcommand{\lra}{\longrightarrow}
\newcommand{\lms}{\longmapsto}

\newcommand{\N}{\mathcal{N}}
\newcommand{\M}{\mathcal{M}}

\usepackage[latin1]{inputenc}

\def\og{\leavevmode\raise.3ex\hbox{$\scriptscriptstyle\langle\!\langle$~}}
\def\fg{\leavevmode\raise.3ex\hbox{~$\!\scriptscriptstyle\,\rangle\!\rangle$}}

\begin{document}

\begin{frontmatter}




%
\selectlanguage{francais}
\title{Feuilletage lisse de $\Sp^5$ par surfaces complexes}

\vspace{-2.6cm} \selectlanguage{english}
\title{Smooth foliation of $\Sp^5$ by complex surfaces}



\author[authorlabel1]{Guillaume Deschamps}
\ead{guillaume.deschamps@univ-brest.fr}

\address[authorlabel1]{UFR de Math\'ematiques, Universit\'e Bretagne Occidentale,
 29200 Brest cedex.}

\begin{abstract}
In 2002 Meersseman-Verjovsky \cite{MV} constructed a smooth,
codimension-one, foliation on 5-sphere by complex surfaces with
two compact leaves. The aim of this note is to improve their
construction in order to give a smooth foliation on 5-sphere by
complex surfaces with only one compact leaf.


\vskip 0.5\baselineskip

\selectlanguage{francais}
\noindent{\bf R\'esum\'e} \vskip 0.5\baselineskip \noindent En
2002 Meersseman-Verjovsky \cite{MV} ont construit un feuilletage
de codimension un de $\Sp^5$ par feuilles complexes, possédant 2
feuilles compactes. Le but de cette note est d'améliorer la
construction afin de munir la sphère de dimension cinq d'un
feuilletage lisse à feuilles complexes avec une seule feuille
compacte.


\end{abstract}
\end{frontmatter}

\selectlanguage{francais}
\section{Introduction}
La note de Novikov \cite{N} parue en 1964, où il esquissait une
démonstration que tout feuilletage lisse de la 3-sphère par
surface possédait une feuille compacte, a eu un impact
considérable. On pouvait alors penser que la construction de
Lawson \cite{L} d'un feuilletage lisse de codimension un  sur
$\Sp^5$ avec une seule feuille compacte, était optimal du point de
vue du nombre de feuilles compactes. Mais on sait, grâce aux
récents travaux de Meigniez \cite{Me}, qu'il existe un feuilletage
lisse de codimension un sur la sphère $\Sp^5$ dont toutes les
feuilles sont denses.

 Trente ans plus tard
Meersseman-Verjovsky \cite{MV} en modifiant la construction de
Lawson ont pu définir un feuilletage lisse à feuilles complexes
sur $\Sp^5$ munit de deux feuilles compactes. Pour les
feuilletages lisses à feuilles complexes on n'a pas de résultat du
type de celui de Meigniez. La question qui se pose alors est donc:
"ce dernier feuilletage est-il optimal en terme de nombre de
feuilles compactes?" Le but de cet article est de montrer que non.

\quad\\
 {\bf Théorème} : {\it Il existe un feuilletage lisse à
feuilles complexes et de codimension un sur $\Sp^5$ ne contenant
qu'une seule feuille compacte.}

\quad\\
Pour démontrer ce théorème nous modifierons légèrement la
construction de \cite{MV} dont nous rappelons ici les notations.

\section{Notations}
On considère $\Sp^5$ comme la sphère unité de $\C^3$. Soit
$$
W=\{z\in\C^3-\{0\}/P(z)=z_1^3+z_2^3+z_3^3=0\}.
$$
C'est une variété complexe. Soit $K$ l'intersection de $W$ avec
$\Sp^5$. On décompose la 5-sphère en deux variétés à bord :
$\mathcal N$, un voisinage tubulaire fermé de $K$ dans $\Sp^5$ et
$\mathcal M$ l'adhérence du complémentaire de $\N$ dans $\Sp^5$.
En particulier le bord commun de $\N$ et de $\M$ est difféomorphe
à $\Sp^1\times K$.

Soit $X$  une variété à bord dont le bord $\partial X$ est une
variété complexe. On rappelle qu'un feuilletage de $X$ par
variétés complexes est dit {\it plat} \cite{MV} s'il s'étend en un
feuilletage à feuilles complexes de $X\cup\partial X\times[0,1]$
qui coïncide avec le feuilletages naturelles du collier $\partial
X\times[0,1]$. L'intérêt de cette définition provient de la
proposition suivante.

\quad\\
{\bf Proposition \cite{MV}} : {\it Soient $(X_i,\mathcal F_i)$
deux feuilletages à bord ($i=1,2$). Supposons les bords
biholomorphes et les feuilletages plats. Alors pour tout
biholomorphisme $\psi$ de $\partial X_1$ sur $\partial X_2$, il
existe un feuilletage par variétés complexes sur l'union
$X_1\cup_\psi X_2$ (recollé le long du bord via $\psi$) dont la
restriction à $X_1$ (respectivement à $X_2$) est $\mathcal F_1$
(respectivement $\mathcal F_2$).}

\quad\\
 On sait que $\N$ admet un feuilletage à feuilles complexes
plat dont la seule feuille compacte est son bord $\partial
\N$
\cite{MV}. La démonstration du théorème revient donc à construire
un feuilletage plat par variétés complexes sur $\M$ avec comme
seule feuille compacte $\partial \M$.

\section{Construction d'un feuilletage plat sur $\M$}
On définit $Y=P^{-1}([0,+\infty[)-\{(0,0,0)\}\subset\C^3$ et
l'application
$$g:(z,t)\in\C^3\times\R\lms P(z)-\phi(t)\in\R$$
 où $\phi$ est la fonction plate en zéro donnée par:
$$
\begin{array}{cccc}
\phi:&]-\infty,1]&\lra&\R\\
&t&\lms& \left\{
\begin{array}{ll}
0&\textrm{ si $t\leq0$}\\
 e^{-\frac{1}{e^-\frac{1}{t}}}&\textrm{ si $t\in
]0,1]$}
\end{array}
\right.
\end{array}
$$
La fonction $\phi$ a les propriétés suivantes:

\begin{enumerate}
\item[i)] $\phi$ est de classe $\mathcal C^\infty$

\item[ii)] $\phi'(t)=\frac{\phi(t)}{t^2e^{-\frac{1}{t}}}>0$ si
$t>0$

\item[iii)] $\phi$ est une bijection de $]0,1]$ sur $]0,\phi(1)]$
d'inverse la fonction
$\phi^{-1}(t)=\frac{1}{\ln\Big(\ln(\frac{1}{t})\Big)}$.
\end{enumerate}
On  prolonge alors $\phi$ sur $[1,+\infty[$ en une fonction
$\mathcal C^\infty$, surjective de $\R$ sur $\R$ et tel que
$\phi'(t)>0\; \forall t>0$.

\quad\\
 Posons alors $\Xi=g^{-1}(\{0\})-\{(0,0,0)\}\times\R$ et
 remarquons que $\Xi$ est l'union de
 $$\Xi^-=g^{-1}(\{0\})\cap\{(z,t)\in\C^3\times\R/t\leq0,
 z\neq[0,0,0)\}
 $$
difféomorphe à $W\times]-\infty,0]$ et de
$$\Xi^+=g^{-1}(\{0\})\cap\{(z,t)\in\C^3\times\R/t\geq0,
 z\neq[0,0,0)\}
 $$
 difféomorphe à $Y$. L'intersection de ces deux pièces est
 difféomorphe à $W=\partial Y$ si bien que $\Xi$ est difféomorphe à $Y$
 augmenté d'un collier infini. On feuillette alors $\Xi^+$ par les niveaux
$$
L_t=\{(z,t)\in\Xi^+/P(z)=\phi(t)\}
$$
et $\Xi^-$ par les niveaux
$$
L_t=\{(z,t)\in\Xi^-/P(z)=\phi(t)=0\}
$$
C'est un feuilletage lisse à feuilles complexes sur $\Xi$. Pour
$0<\lambda<1$, on note
$$
\begin{array}{cccc}
G:&\Xi&\lra&\Xi\\
&(z,t)&\lms&\Big(\lambda jz, h(t)\Big)
\end{array}
$$
où $h:\R\lra \R$ est la fonction constante égale $t$ sur
$]-\infty,0]$ et sur $\R^+$ égale à:
$$h^+(t)=\phi^{-1}\Big(\lambda^3\phi(t)\Big).$$
Le groupe engendré par $G$ agit librement, proprement sur $\Xi$,
respecte le feuilletage et est holomorphe en restriction aux
feuilles. Le quotient $Y_1$ est donc une variété feuilletée par
feuilles complexes.

 \quad\\
{\bf Lemme 1.} : {\it La variété  $Y_1$ est difféomorphe à
$\M\cup\partial\M\times]-\infty,0]$.}

\quad\\
{\it Preuve.}  Si on note $int(\Xi^+)$ l'intérieur de $\Xi^+$
alors le difféomorphisme :
$$
\begin{array}{ccc}
int(\Xi^+)&\lra&P^{-1}(1)\times]0,+\infty[\\
(z,t)&\lms&\Big(\frac{z}{\phi^\frac{1}{3}(t)},\phi(t)\Big)
\end{array}
$$
induit un difféomorphisme entre les feuilletages naturels de ces
deux variétés. De plus il conjugue $G$ à :
$$
\begin{array}{lccc}
\tilde G: &P^{-1}(1)\times]0,+\infty[&\lra&P^{-1}(1)\times]0,+\infty[\\
&(z,t)&\lms&(jz,\lambda^3t)
\end{array}
$$
Le quotient de  $int(\Xi^+)$ par $G$ est donc difféomorphe à un
fibré en cercle de fibre $P^{-1}(1)$ et de monodromie donnée par
la multiplication par $j$.  Maintenant la fibration de Milnor qui
envoie un point $z$ de l'intérieur de $\M$ sur $P(z)/\vert
P(z)\vert$ a la même monodromie \cite{M}. On a bien $int(\Xi^+)/G$
difféomorphe à $int(\M)$ et donc $Y_1$ difféomorphe à
$\M\cup\partial\M\times]-\infty,0]$. $\square$

\quad\\
 Le feuilletage que nous avons construit sur $Y_1$ est
lisse du fait du choix de la fonction $\phi$:

\quad\\
{\bf Lemme 2.} : {\it La fonction $h$ est de classe $\mathcal
C^\infty$ en zéro.}

\quad\\
 Ce lemme nous dit précisément que nous
avons construit un feuilletage
 plat à feuilles complexes sur $\M$ dont la seule feuille compacte
 est son bord. Ce qui
conclut la démonstration de notre théorème.

\quad\\
{\it Preuve du Lemme 2.} On pose
$u(t)=1-\ln(\lambda^3)e^{-\frac{1}{t}}$, on peut alors écrire
$$\forall t\in]0,1]\;h^+(t)=\phi^{-1}\Big(\lambda^3\phi(t)\Big)=\frac{t}{t\ln\Big(1-\ln(\lambda^3)e^{-\frac{1}{t}}\Big)+1}=
\frac{t}{t\ln(u(t))+1}
$$
Mais en zéro on a $e^{-\frac{1}{t}}=o(t^n),\;\forall
n\in\mathbb{N}$ de sorte que $u(t)=1+o(t^n),\;\forall
n\in\mathbb{N}$ et donc:
$$
 h^+(t)=\frac{t}{t\ln(1+o(t^n))+1}
 =\frac{t}{o(t^{n+1})+1}
 =t+o(t^{n+1})
 $$
En d'autres termes on a $h^+(0)=0,\; {h^+}'(0)=1$ et
${h^+}^{(n)}(0)=0,\;\forall n>1$. La fonction $h$ est bien de
classe $\mathcal C^\infty$ en zéro. $\square$

\quad\\
 {\bf Remerciements.} Je tiens à remercier L. Meersseman d'avoir
 porté mon attention sur cette question.
  Ce travail a bénéficié d'une aide de l'Agence Nationale
 de la Recherche portant la référence ANR-08-JCJC-0130-01.

\end{document}